\documentclass[12pt]{article}
\usepackage{amssymb}
\usepackage{amsmath}
\usepackage{gdgspace}

\numberwithin{equation}{section}

\begin{document}
\thispagestyle{empty}
\begin{center}
\textbf{TRANSFORMATIONS OF RAMANUJAN'S SUMMATION FORMULA AND ITS APPLICATIONS}
\end{center}
\vspace{.25cm}
\begin{center}
Chandrashekar Adiga\footnote{$\text{e-mail} : \text{c}_{-}\text{adiga@hotmail.com}$} and N.Anitha \footnote{$\text{e-mail} : \text{ani}_{-}\text{nie@indiatimes.com}$}\\
Department of Studies in Mathematics\\
University of Mysore, Manasagangotri\\
Mysore-570006, INDIA.\\
\end{center}
\begin{center}
Taekyun Kim \footnote{$\text{e-mail} : \text{tkim@kongju.ac.kr}$}\\
Institute of Science Education,\\
Kongju National University,\\
Kongju 314-701,\\
KOREA.
\end{center}
\vspace{.25cm}
\textbf{Abstract:}\quad In this paper, we obtain some new transformation formulas for Ramanujan's $_1\psi_1$ summation formula and also establish some eta-function identities. We also deduce a $q$- gamma function identity, an $q$-integral and some interesting series representations for $\frac{\pi^{3/2}}{2\sqrt{2}\Gamma^2{(3/4)}}$ and the beta function $B(x,y)$ .
\vspace{.25cm}

\begin{tabbing}
                               \> $q$-gamma function, $q$-integrals.
\end{tabbing}
\vspace{.25cm}
\begin{tabbing}
\end{tabbing}
\vspace{.5cm}

\section{Introduction :}
One of the most celebrated identities of Ramanujan is his $_1\psi_1$-summation formula which can
be stated as follows:
\begin{equation}
\sum_{n=-\infty}^{\infty}\frac{(a)_n}{(b)_n}z^n = \frac{(az)_\infty(q/az)_\infty(q)_\infty(b/a)_\infty}{(z)_\infty(b/az)_\infty(b)_\infty(q/a)_\infty}
\end{equation}
where $\mid b/a \mid < \mid z \mid < 1.$  This identity was first brought before the mathematical world by G. H. Hardy [7] who described it as ``A remarkable formula with many parameters". Hardy did not supply a proof but indicated that a proof could be constructed from\linebreak the $q$-binomial theorem.  There are several proofs of (1.1) in the literature [c.f 1 ]. Most proofs of (1.1) found in the literature do depend on the $q$-binomial theorem. The transformations of the basic hypergeometric series are extremely useful in the theory of partitions. The main purpose of this paper is to obtain some transformations for Ramanujan's $_1\psi_1$ summation formula (1.1).

\qquad As customary we employ the following definitions and notations:
\begin{equation}
 (a)_\infty:= (a;q)_\infty = \prod_{n=0}^{\infty}(1-aq^n), \quad \mid q \mid < 1,
\end{equation}
\begin{equation}
(a)_n:= \frac{(a)_\infty}{(aq^n)_\infty},\quad n: \text{any integer}.
\end{equation}
The generalized basic hypergeometric series is defined by
\begin{equation}
_r\phi_s\left[\begin{array}{ll} a_1,a_2,...,a_r;q;z\\ b_1,b_2,...,b_s
                 \end{array}
  \right] = \sum_{n=0}^{\infty}\frac{(a_1)_n(a_2)_n...(a_r)_n}{(q)_n(b_1)_n(b_2)_n...(b_s)_n}
z^n,\quad  \mid z \mid < 1, \quad \mid q \mid < 1.
\end{equation}\\

The bilateral basic hypergeometric series is defined by
\begin{equation}
_r\psi_r\left[\begin{array}{ll} a_1,a_2,...,a_r;q;z\\ b_1,b_2,...,b_r
                 \end{array}
  \right] = \sum_{n=-\infty}^{\infty}\frac{(a_1)_n(a_2)_n...(a_r)_n}{(b_1)_n(b_2)_n...(b_r)_n}
z^n,
\end{equation}
where $\left|\frac{b_1...b_r}{a_1...a_r}\right|$ $<$ $\mid z \mid$ $<$ 1, $\mid q \mid$
$<$ 1.\\

In Section 3 we obtain some special cases of the transformations for Ramanujan's $_1\psi_1$
summation formula obtained in Section 2 and in Section 4 we also deduce some eta-function identities.
In the final section we derive some interesting series representations for the beta function $B(x,y)$ and
$\frac{\pi^{3/2}}{2\sqrt2\Gamma^2{(3/4)}}$.

\section{Main Results}
To prove our main result we need the Heine's transformation[5, Eq.(III.1) p.241, Eq.(III.2),
p.241]
\begin{eqnarray}
\sum_{n=0}^{\infty}\frac{(a)_n(b)_n}{(q)_n(c)_n}z^n &=& \frac{(b)_\infty(az)_\infty}{(c)_\infty
(z)_\infty} \sum_{n=0}^{\infty}\frac{(c/b)_n(z)_n}{(q)_n(az)_n}b^n\\
& =&\frac{(c/b)_\infty(bz)_\infty}{(c)_\infty(z)_\infty}\sum_{n=0}^{\infty}
\frac{(abz/c)_n(b)_n}{(q)_n(bz)_n}(c/b)^n.
\end{eqnarray}

\textbf{Theorem 2.1.}
If  $\mid b/a \mid  <  \mid z \mid  <  1$ and   $\mid q \mid <$   min $\{1, \quad \mid b \mid \}$
 \quad then,
\begin{eqnarray}
\sum_{n=-\infty}^{\infty}\frac{(a)_n}{(b)_n}z^n  = -1 +\frac{(b/a)_\infty(az)_\infty}{(b)_\infty(z)_\infty}\sum_{n=0}^{\infty}\frac{(a)_n(aqz/b)_n}{(q)_n(az)_n}(b/a)^n  \nonumber\\
 +\frac{(q/b)_\infty(bq/az)_\infty}{(q/a)_\infty(b/az)_\infty} \sum_{n=0}^{\infty}\frac{(b/a)_n(b/az)_n}{(q)_n(bq/az)_n}(q/b)^n.
\end{eqnarray}
\textbf{Proof:} We have
\begin{eqnarray}
\sum_{n=-\infty}^{\infty}\frac{(a)_n}{(b)_n}z^n & = &\sum_{n=0}^{\infty}\frac{(a)_n}{(b)_n}z^n                               +\sum_{n=1}^{\infty}\frac{(q/b)_n}{(q/a)_n}(b/az)^n \nonumber \\
                                                 & = &-1 + \sum_{n=0}^{\infty}\frac{(a)_n}{(b)_n}z^n
+\sum_{n=0}^{\infty}\frac{(q/b)_n}{(q/a)_n}(b/az)^n.
\end{eqnarray}
Putting $a=q$, $b=a$ and $c=b$ in (2.2) we obtain
\begin{equation}
\sum_{n=0}^{\infty}\frac{(a)_n}{(b)_n}z^n = \frac{(b/a)_\infty(az)_\infty}{(b)_\infty(z)_\infty
}\sum_{n=0}^{\infty}\frac{(a)_n(aqz/b)_n}{(q)_n(az)_n}(b/a)^n.
\end{equation}
Changing $a$ to $q$, $b$ to $q/b$, $c$ to $q/a$ and $z$ to $b/az$ in (2.1) we obtain
\begin{equation}
\sum_{n=0}^{\infty}\frac{(q/b)_n}{(q/a)_n}(b/az)^n =
\frac{(q/b)_\infty(bq/az)_\infty}{(q/a)_\infty(b/az)_\infty} \sum_{n=0}^{\infty}\frac{(b/a)_n(b/az)_n}{(q)_n(bq/az)_n}(q/b)^n.
\end{equation}\\
Substituting (2.5) and (2.6) in (2.4) we obtain (2.3).\\\\

\textbf{Theorem 2.2.}
If $\mid \frac{b}{a}\mid < \mid z\mid < 1$ and $\mid a\mid < 1$
then,
\begin{eqnarray}
\sum_{n=-\infty}^{\infty}\frac{(a)_n}{(b)_n}z^n  = -1 +\frac{(a)_\infty(qz)_\infty}{(b)_\infty(z)_\infty}\sum_{n=0}^{\infty}\frac{(b/a)_n(z)_n}{(q)_n(qz)_n}(a)^n  \nonumber\\
 +\frac{(b/a)_\infty(q/az)_\infty}{(q/a)_\infty(b/az)_\infty}\sum_{n=0}^{\infty}\frac{(q/z)_n(q/b)_n}{(q)_n(q/az)_n}(b/a)^n.
\end{eqnarray}
\textbf{Proof:}
 Changing $a$ to $q$, $b$ to $a$ and then $c$ to $b$ in (2.1), we deduce
\begin{equation}
\sum_{n=0}^{\infty}\frac{(a)_n}{(b)_n} z^n = \frac{(a)_\infty(qz)_\infty}{(b)_\infty(z)_\infty}\sum_{n=0}^{\infty}\frac{(b/a)_n(z)_n}{(q)_n(qz)_n}(a)^n .
\end{equation}
Putting $a=q$, $b=q/b$, $c=q/a$ and $z=b/az$ in (2.2), we obtain
\begin{equation}
\sum_{n=0}^{\infty}\frac{(q/b)_n}{(q/a)_n}(b/az)^n = \frac{(b/a)_\infty(q/az)_\infty}{(q/a)_\infty(b/az)_\infty}\sum_{n=0}^{\infty}\frac{(q/z)_n(q/b)_n}{(q)_n(q/az)_n}(b/a)^n.
\end{equation}
On employing (2.8) and (2.9) in (2.4) we obtain Theorem 2.2.\\

\textbf{Theorem 2.3.}
If $\mid \frac{b}{a}\mid < \mid z\mid < 1$ then,
\begin{eqnarray}
\sum_{n=-\infty}^{\infty}\frac{(a)_n}{(b)_n}z^n  = -1 +\frac{(b/a)_\infty(az)_\infty}{(b)_\infty(z)_\infty}\sum_{n=0}^{\infty}\frac{(aqz/b)_n(a)_n}{(q)_n(az)_n}(b/a)^n  \nonumber\\
 +\frac{(b/a)_\infty(q/az)_\infty}{(q/a)_\infty(b/az)_\infty}\sum_{n=0}^{\infty}\frac{(q/z)_n(q/b)_n}{(q)_n(q/az)_n}(b/a)^n.
\end{eqnarray}
\textbf{Proof:}
On employing (2.5) and (2.9) in (2.4) we obtain Theorem 2.3.\\
\section{Special Cases}

Putting $a=-1/q$ and $b=-1$ in (2.3), we obtain
\begin{equation}
\sum_{n=-\infty}^{\infty}\frac{z^n}{1+q^{n-1}} = \frac{-q}{1+q} + \frac{q}{1+q}\frac{(q)_\infty(-z/q)_\infty}{(-1)_\infty(z)_\infty}\sum_{n=0}^{\infty}\frac{(z)_n(-1/q)_n}{(q)_n(-z/q)_n}q^n + q \sum_{n=0}^{\infty}\frac{1}{1-\frac{q^{n+1}}{z}}(-q)^n,\\
\end{equation}\\
\begin{text} where  $\mid q \mid < \mid z \mid < 1 $ \end{text}.\\

Changing $b$ to $-q^3$, $a$ to $-q$,  $z$ to $q$ in (2.7) we get,
\begin{eqnarray}
\sum_{n=-\infty}^{\infty}\frac{q^{n}}{(1+q^{n+1})(1+q^{n+2})} = \frac{(1+q^2)(1+q)}{q(1-q)}.
\end{eqnarray}\\

Putting $b=-q^2$, $a=-1/q$ and $z=q^{1/2}$ in  (2.10) and then changing $q$ to $q^2$, we obtain
\begin{eqnarray}
\sum_{n=-\infty}^{\infty}\frac{2(1+1/q^2)(1+q^2)q^n}{(1+q^{2n-2})(1+q^{2n})(1+q^{2n+2})}=-1+\frac{(q^6;q^2)_\infty(-1/q;q^2)_\infty}{(-q^4;q^2)_\infty(q;q^2)_\infty}\nonumber\\
\sum_{n=0}^{\infty}\frac{(1/q^3;q^2)_n(-1/q^2;q^2)_n}{(q^2;q^2)_n(-1/q;q^2)_n}q^{6n}+\frac{(q^6;q^2)_\infty(-q^3;q^2)_\infty}{(-q^4;q^2)_\infty(q^5;q^2)_\infty}\sum_{n=0}^{\infty}\frac{(q;q^2)_n(-1/q^2;q^2)_n}{(q^2;q^2)_n(-q^3;q^2)_n}q^{6n}.
\end{eqnarray}\\
\pagebreak
\section{Eta-function identities}

The Dedekind eta-function is defined by

$$\eta(\tau) :=  e^{\pi i \tau/12}\prod_{n=1}^{\infty}(1-e^{2\pi i n \tau})$$
\begin{equation}
\qquad\qquad  = q^{1/24}(q;q)_\infty,
\end{equation}
where $q = e^{2\pi i \tau}$ and Im $\tau$ $>$ 0.\\

The eta-function is useful in the study of modular forms. Berndt and L. C. Zhang [3] have
obtained a number of eta-function identities found in
Ramanujan's notebook employing the theory of modular forms.
In [4], Bhargava and Somashekara show how a family of interesting eta-function
identities can be obtained from a limiting case of the Ramanujan's $_1\psi_1$-summation formula.

Now we derive some new $\eta$- function identities.\\

Employing Ramanujan's $_1\psi_1$ summation formula in (2.3) and then
putting $a=1/q$, $b=q^{1/2}$, $z=q^{1/2}$ and changing $q$ to $q^2$, we deduce\\

\begin{eqnarray}
\frac{\eta{(\tau)}}{\eta^2(2\tau)}= q^{-1/8}-\frac{q^{7/8}}{1+q}\frac{(q;q^2)_\infty}{(q^2;q^2)_\infty}\sum_{n=0}^{\infty}\frac{(q^3;q^2)_n}{(q^4;q^2)_n}q^n.
\end{eqnarray}\\

This can also be obtained from $q$-binomial theorem.\\

Employing Ramanujan's $_1\psi_1$ summation formula in (2.7) and then
putting $a=-q$, $b=-q^3$, $z=q^{1/2}$ and changing $q$ to $q^2$, we deduce\\

\begin{eqnarray}
\frac{\eta^{10}(2\tau)}{\eta^4(\tau)\eta^2(4\tau)}= \frac{2(1+q)q^{4/3}}{(1+q^2)(1+q^4)}-\frac{2q^{4/3}}{1-q}\sum_{n=0}^{\infty}\frac{1-q^{2n+2}}{1-q^{2n+1}}(-q^2)^n-\nonumber\\
\frac{2(1+q)q^{4/3}}{(1+q^2)(1+q^4)}\frac{(q^4;q^2)_\infty(-1/q;q^2)_\infty}{(-1;q^2)_\infty(q^3;q^2)_\infty}\sum_{n=0}^{\infty}\frac{(q;q^2)_n(-1/q^4;q^2)_n}{(q^2;q^2)_n(-1/q;q^2)_n}q^{4n}.
\end{eqnarray}\\

Employing Ramanujan's $_1\psi_1$ summation formula in (2.3) and then
putting $a=-q^{-5/3}$, $b=-q^{1/3}$, $z=q^{4/3}$ and changing $q$ to $q^3$, we deduce\\

\begin{eqnarray}
\frac{\eta^3(3\tau)}{\eta(\tau)}=\frac{q^{4/3}(1+q+q^2)}{(1+q^2)(1+q^5)}\left[-1+\frac{(q^6;q^3)_\infty(-1/q;q^3)_\infty}{(-q;q^3)_\infty(q^4;q^3)_\infty}\sum_{n=0}^{\infty}\frac{(q;q^3)_n(-1/q^5;q^3)_n}{(q^3;q^3)_n(-1/q;q^3)_n}q^{6n}\right]\nonumber\\
+\frac{q^{4/3}(1+q+q^2)}{(1+q^2)(1+q^5)}\frac{(q^6;q^3)_\infty(-q^4;q^3)_\infty}{(-q^8;q^3)_\infty(q^2;q^3)_\infty}\sum_{n=0}^{\infty}
\frac{(1/q;q^3)_n(-q^2;q^3)_n}{(q^3;q^3)_n(-q^4;q^3)_n}q^{6n}.
\end{eqnarray}

\section{Some q-Gamma Function Identities:}
\qquad F.H. Jackson [8] defined the $q$-analogue of the gamma function by
\begin{equation}
\Gamma_q(x)=\frac{(q)_\infty}{(q^x)_\infty}(1-q)^{1-x}, \quad  0 < q < 1.
\end{equation}
Jackson [8] also defined a $q$-integral by,
\begin{equation}
\int_{0}^{a} f(t)d_q(t) = a(1-q)\sum_{n=0}^{\infty} f(aq^n) q^n
\end{equation}
and
\begin{equation}
\int_{0}^{\infty}f(t)d_q(t) = (1-q)\sum_{n=-\infty}^{\infty}f(q^n).q^n .
\end{equation}\\

\qquad In this section, as an application of the transformation formulas of Section 2, we deduce some $q$-gamma function identities and an $q$-integral which give interesting series representations for the beta function B(x,y) and $\frac{\pi^{3/2}}{2\sqrt2\Gamma^2(3/4)}$.\\

\textbf{Theorem 5.1.}
If $0 < z < b-a <1$ and $b<1$ then,
\begin{eqnarray}
\frac{\Gamma_q(b)\Gamma_q(1-a)\Gamma_q(z)\Gamma_q(b-a-z)}{\Gamma_q(b-a)\Gamma_q(a+z)\Gamma_q(1-a-z)}=-(1-q)^{a+1-b}+(1-q)^{a+1-b}\frac{\Gamma_q(b)\Gamma_q(z)}{\Gamma_q(b-a)\Gamma_q(a+z)}\nonumber\\
\sum_{n=0}^{\infty}\frac{(q^{a+1+z-b})_n(q^a)_n}{(q)_n(q^{a+z})_n}(q^{b-a})^n
+ \frac{\Gamma_q(1-a)\Gamma_q(b-a-z)}{\Gamma_q(1-b)\Gamma_q(b+1-a-z)}
\sum_{n=0}^{\infty}\frac{(q^{b-a})_n(q^{b-z-a})_n}{(q^{b+1-a-z})_n(q)_n}(q^{1-b})^n.\nonumber\\
\end{eqnarray}
\textbf{Proof:}\qquad
Changing $a$ to $q^a$, $b$ to $q^b$ and $z$ to $q^z$ in (2.3) we obtain the required result.\\

\textbf{Corollary:}\qquad
In Theorem 5.1 if we let $q \rightarrow 1$ and put $a = 0$ we deduce,
\begin{eqnarray}
\frac{\Gamma(1-b)\Gamma(b+1-z)}{\Gamma(1-z)} =\sum_{n=0}^{\infty}\frac{(b)_n}{n!}\frac{b-z}{b-z+n}.
\end{eqnarray}
\textbf{Remark:}\qquad Putting $b-z=y$ and $1-b=x$ in (5.5) we obtain a well known [6]
representation for the beta function $B(x,y)$ given by,
\begin{equation}
B(x,y) = \sum_{n=0}^{\infty}\frac{\prod_{k=1}^{n}(k-x)}{n!}\frac{1}{n+y}.
\end{equation}

\textbf{Theorem 5.2.} If $0 < z < b-a<1$ and $a >0 $  then,
 \begin{eqnarray}
\frac{\Gamma(a)\Gamma(z)\Gamma(1-a)\Gamma(b-a-z)}{\Gamma(a+z)\Gamma(b-a)\Gamma(1-a-z)} = \sum_{n=0}^{\infty}\frac{(b-a)_n}{ n!(n+z)}.
\end{eqnarray}

\textbf{Proof:}\qquad Changing $a$ to $q^a$, $b$ to $q^b$ and $z$ to $q^z$ in (2.7)
we obtain
\begin{eqnarray}
\frac{\Gamma_q(z)\Gamma_q(b-a-z)\Gamma_q(b)\Gamma_q(1-a)}{\Gamma_q(a+z)\Gamma_q(1-a-z)\Gamma_q(b-a)}=-(1-q)^{a+1-b} + \frac{\Gamma_q(b)\Gamma_q(z)}{\Gamma_q(a)\Gamma_q(z+1)}\nonumber\\ \sum_{n=0}^{\infty}\frac{(q^{b-a})_n(z)_n}{(q)_n(q^{1+z})_n}(q^{a})^n +(1-q)^{a+1-b}\frac{\Gamma_q(1-a)\Gamma_q(b-a-z)}{\Gamma_q(1-a-z)\Gamma_q(b-a)}\sum_{n=0}^{\infty}\frac{(q^{1-z})_n(q^{1-b})_n}{(q)_n(q^{1-a-z})_n}(q^{b-a})^n.\nonumber\\
\end{eqnarray}

Letting $q \rightarrow 1$ in (5.8) we obtain Theorem 5.2.\\

\textbf{Remark:}
Putting $a=1/4$, $b=1$ and $z=1/2$ in Theorem 5.2 we obtain,
\begin{eqnarray}
\frac{\pi^{3/2}}{2\sqrt2\Gamma^2{(3/4)}} = \sum_{n=0}^{\infty}\frac{(1/2)_n}{n!(4n+1)}
\end{eqnarray}
which is very similar to the result obtained by Ramanujan[2, p.24].\\

\textbf{Theorem 5.3.}
If $0 < z < b-a <1$  then,

\begin{eqnarray}
\frac{\Gamma_q(b)\Gamma_q(1-a)\Gamma_q(z)\Gamma_q(b-a-z)}{\Gamma_q(b-a)\Gamma_q(a+z)\Gamma_q(1-a-z)}&=&-(1-q)^{a+1-b}+\frac{\Gamma_q(1-a)\Gamma_q(b-a-z)}{\Gamma_q(b-a)\Gamma_q(1-z)\Gamma_q(1-b)}\nonumber\\
\int_{0}^{1}f(x)d_q(x) &+ &\frac{\Gamma_q(b)\Gamma_q(z)}{\Gamma_q(b-a)\Gamma_q(a+1+z-b)\Gamma_q(a)}\int_{0}^{1}g(x)d_q(x)\nonumber\\
\end{eqnarray}
where,
$$f(x)= \frac{(xq)_\infty(xq^{1-a-z})_\infty}{(xq^{1-z})_\infty(xq^{1-b})_\infty}x^{b-a-1}$$ and
$$g(x) = \frac{(xq)_\infty(xq^{a+z})_\infty}{(xq^{a+1+z-b})_\infty(xq^a)_\infty} x^{b-a-1}.$$ \\

\textbf{Proof:}\qquad
Changing $a$ to $q^a$, $b$ to $q^b$ and $z$ to $q^z$ in (2.10) we get,
\begin{eqnarray}
\frac{\Gamma_q(b)\Gamma_q(1-a)\Gamma_q(z)\Gamma_q(b-a-z)}{\Gamma_q(b-a)\Gamma_q(a+z)\Gamma_q(1-a-z)}(1-q)^{b-a-1}= -1+ \frac{\Gamma_q(1-a)\Gamma_q(b-a-z)}{\Gamma_q(b-a)\Gamma_q(1-z)\Gamma_q(1-b)}\nonumber\\(1-q)^{(b-a)}
\sum_{n=0}^{\infty}\frac{(q^{n+1})_\infty(q^{n+1-a-z})_\infty}{(q^{n+1-z})_\infty(q^{n+1-b})_\infty}(q^{b-a})^{n}+\frac{\Gamma_q(b)\Gamma_q(z)}{\Gamma_q(b-a)\Gamma_q(a+1+z-b)\Gamma_q(a)}\nonumber\\
(1-q)^{(b-a)}\sum_{n=0}^{\infty}\frac{(q^{n+1})_\infty(q^{n+a+z})_\infty}{(q^{n+a+1+z-b})_\infty(q^{n+a})_\infty}(q^{b-a})^n\nonumber
\end{eqnarray}

\begin{equation}
\end{equation}
Applying (5.2) to the right side of (5.11) we obtain Theorem 5.3.

\textbf{Corollary:}\qquad
In Theorem 5.3 if we let $q \rightarrow 1$ we deduce,
\begin{eqnarray}
\frac{\Gamma(b)\Gamma(1-a)\Gamma(z)\Gamma(b-a-z)}{\Gamma(a+z)\Gamma(1-a-z)} =\nonumber\\
\left[\frac{\Gamma(1-a)\Gamma(b-a-z)}{\Gamma(1-z)\Gamma(1-b)}+
\frac{\Gamma(b)\Gamma(z)}{\Gamma(a)\Gamma(a+1+z-b)}\right] B(b-a,a-b+1).\nonumber\\
\end{eqnarray}
\begin{center}

\end{center}

\begin{thebibliography}{99}
\bibitem{}
C. Adiga, B. C. Berndt, S. Bhargava and G. N. Watson, Chapter 16 of Ramanujan's \linebreak
second notebook: {\em Theta-functions and q-series}, Mem.Amer.Math.Soc.No.315, 53(1985), 1-85.
\bibitem{}
B. C. Berndt, {\em Ramanujan's Notebooks II}, Springer-Verlag, New York, 1991.
\bibitem{}
B. C. Berndt and L. C. Zhang, {\em Ramanujan's identities for eta functions}, Math.Ann.292(3)
(1992), 561-573.
\bibitem{}
S. Bhargava and D. D. Somashekara, {\em Some eta function identities deducible from \linebreak Ramanujan's $_1\psi_1$Summation}, J.Math.Anal.Appl.176(2)(1993), 554-560.
\bibitem{}
G. Gasper and M. Rahman, {\em Basic Hypergeometric Series}, Cambridge University Press, 1990.
\bibitem{}
I. S. Gradshteyn and I. M. Ryzhik, {\em Tables of Integrals}, Series and Products, 4th ed.,
Academic Press, New York, 1965.
\bibitem{}
G. H. Hardy, {\em Ramanujan}, 3rd ed. Chelsea, NewYork, 1978.
\bibitem{}
F. H. Jackson, {\em On q-definite integrals}, Quar.Jour.Pure and Applied Math. 41(1910), 193-203.
\end{thebibliography}
\end{document}